\documentclass{amsart}

\usepackage{hyperref}
\hypersetup{
colorlinks=true,
linkcolor=blue,
filecolor=blue,      
urlcolor=blue,
citecolor=blue,
}

\usepackage{xypic}
\input xy
\xyoption{all}
\usepackage{epsfig}
\usepackage{amsthm}
\usepackage{amssymb}
\usepackage{amsmath}
\usepackage{amscd}
\usepackage{color}

\newcommand{\lb}{\varLambda}

\newcommand{\A}{\mathcal{A}}
\def \<{\langle}
\def \>{\rangle}

\newcommand{\bg}{\begin{equation}}
\newcommand{\ed}{\end{equation}}
\newcommand{\bga}{\begin{eqnarray}}
\newcommand{\eda}{\end{eqnarray}}

\def\cbdu{\par{\raggedleft$\Box$\par}}

\newtheorem {Theorem}  {Theorem}

\numberwithin{Theorem}{section}

\newtheorem {Lemma}[Theorem]  {Lemma}
\newtheorem {Conjecture}[Theorem]{Conjecture}
\newtheorem {Question}[Theorem]{Question}
\theoremstyle{definition}

\theoremstyle{remark}

\def \l{\lambda}
\expandafter\chardef\csname pre amssym.def
at\endcsname=\the\catcode`\@ \catcode`\@=11
\def\undefine#1{\let#1\undefined}
\def\newsymbol#1#2#3#4#5{\let\next@\relax
 \ifnum#2=\@ne\let\next@\msafam@\else
 \ifnum#2=\tw@\let\next@\msbfam@\fi\fi
 \mathchardef#1="#3\next@#4#5}
\def\mathhexbox@#1#2#3{\relax
 \ifmmode\mathpalette{}{\m@th\mathchar"#1#2#3}%
 \else\leavevmode\hbox{$\m@th\mathchar"#1#2#3$}\fi}
\def\hexnumber@#1{\ifcase#1 0\or 1\or 2\or 3\or 4\or 5\or 6\or 7\or 8\or
 9\or A\or B\or C\or D\or E\or F\fi}

\font\teneufm=eufm10 \font\seveneufm=eufm7 \font\fiveeufm=eufm5
\newfam\eufmfam
\textfont\eufmfam=\teneufm \scriptfont\eufmfam=\seveneufm
\scriptscriptfont\eufmfam=\fiveeufm

\catcode`\@=\csname pre amssym.def at\endcsname

\newcounter{remark}
\def \grad {\nabla}

\DeclareMathOperator{\curl}{curl}

\renewcommand{\l}{\lambda}

\newcommand{\s}{\sigma}

\newcommand{\R}{\mathbf{R}}

\def  \R   {{\mathbb R}}
\def  \Z   {{\mathbb Z}}

\def  \T   {{\mathbb T}}

\def  \12  {{\frac{1}{2}}}
\def  \p   {\partial}

\def  \o  {\omega}

\def\build#1_#2^#3{\mathrel{\mathop{\kern 0pt#1}\limits_{#2}^{#3}}}

\begin{document}

\title[The number of degrees of freedom for the 2D NSE]{The number of degrees of freedom for the 2D Navier-Stokes equation: a connection with Kraichnan's theory of turbulence}

\author [Alexey Cheskidov]{Alexey Cheskidov}
\address{Department of Mathematics, Stat. and Comp.Sci.,  University of Illinois Chicago, Chicago, IL 60607,USA}
\address{School of Mathematics, Institute for Advanced Study, Princeton, NJ 08540,USA}
\email{acheskid@uic.edu} 
\author [Mimi Dai]{Mimi Dai}
\address{Department of Mathematics, Stat. and Comp.Sci.,  University of Illinois Chicago, Chicago, IL 60607,USA}
\address{School of Mathematics, Institute for Advanced Study, Princeton, NJ 08540,USA}
\email{mdai@uic.edu}






\begin{abstract}
We estimate the number of degrees of freedom of solutions of the 2D Navier-Stokes equation, proving that its mathematical analog, the number of determining modes, is bounded by the Kraichnan number squared. In particular, this provides new bounds on the number of determining modes in term of the Grashof number for solutions that are not highly intermittent.

\bigskip

KEY WORDS: Navier-Stokes equations, determining modes, Kraichnan's number, global attractor.

\hspace{0.02cm}CLASSIFICATION CODE: 35Q35, 37L30.
\end{abstract}

\maketitle

\section{Introduction}

We consider the Navier-Stokes equation (NSE) on a 2D torus $\mathbb T^2=[0,L]^2$
\begin{equation}
  \label{nse}
  \left\{
    \begin{array}{l}
      u_t + (u \cdot \grad) u - \nu\Delta u + \grad p = f \\
      \grad \cdot u = 0,
    \end{array}
  \right.
\end{equation}
with an external forcing $f$. The unknowns are $u$, the vector velocity field, and $p$, the scalar pressure function. We consider finite energy initial data, and the forcing
$f$ is assumed to have zero mean.  The vorticity $\omega=\nabla\times u=\curl u$ satisfies the equation
\begin{equation}\label{Vnse}
\omega_t+u\cdot\nabla \omega=\nu\Delta \omega + \curl f.
\end{equation}

The purpose of this paper is to estimate the number of degrees of freedom of solutions of the 2D NSE. More precisely, we prove that its mathematical analog, the number of determining modes, is bounded by the Kraichnan number squared $\kappa_\eta^2$, consistent with Kraichnan's theory of two-dimensional turbulence \cite{Kr}.

The notion of determining modes was introduced by  Foias and Prodi in the seminal work \cite{FP} where it was shown that high modes of
a solution to the 2D NSE are controlled by low modes asymptotically as time goes to infinity.  Then the number of such
determining modes $N$ was estimated by Foias, Manley, Temam, and Treve \cite{FMTT} and later improved by Jones and Titi  \cite{JT92,JT}.
We refer the readers to  \cite{CFMT,CFT,FJKT,FJKTgeneral,FJKT12, FMRT,FT,FT84, FT-attractor,FTiti,OT03, OT08} and references therein for more background and related results.

In the remarkable work of Jones and Titi  \cite{JT}, the number of determining modes $N$ was estimated as
\begin{equation} \label{e:JT}
N \lesssim G,
\end{equation}
where $G$ is the Grashof number that measures the size of the force $f$ in terms of the $L^2$ norm, see Section~\ref{s:Long_time_behaviour}. On the other hand, the Hausdorff dimension of the  global attractor $\mathcal A$ for the 2D NSE, which is another way to measure the number of degrees of freedom, was estimated in \cite{CFT} as
 \begin{equation} \label{e:A-dim}
 \mathrm{dim_H}\mathcal A \lesssim G^\frac{2}{3}(1+\log G).
 \end{equation}
Remarkably this bound was proved to be optimal up to a logarithmic correction \cite{L93}, which suggests that the number of degrees of freedom of two-dimensional turbulent flows is of order $G^\frac{2}{3}$. This leads to the following
\begin{Conjecture} \label{c:G-conj}
For two-dimensional turbulent flows, the number of determining modes is bounded by $G^{\frac{2}{3}}$.
\end{Conjecture}

Heuristically, as suggested by the Landau and Lifschitz description of the number of degrees of freedom of a turbulent flow, the modes in the inertial range should determine the behavior of solution of the NSE. In the three-dimensional case, this indicates that Kolmogorov's dissipation wavenumber cubed $\kappa_{\mathrm{d}}^3$, defined in terms of the average energy dissipation rate, should be an appropriate bound for the number of determining modes, which was proved in the time averaged sense employing the wavenumber splitting framework developed in our previous works on the 3D NSE \cite{CD-nse, CDK}. Even though the bound on the time-average of the number of determining modes validates Kolmogorov's predictions, no pointwise bound on the number of determining modes has been obtained yet. This is partially due to the fact that the solutions of the 3D NSE can potentially blow up, but there are other serious obstacles as well.

In the two-dimensional case, Kraichnan's theory of turbulence suggests that the wavenumber $\kappa_\eta$ (Kraichnan's number), defined in terms of the average enstrophy dissipation rate  $\eta$  (see Subsection~\ref{sub:Kraichnan's_number}), indicates the end of the inertial range where the  enstrophy cascade occurs. Hence, it is natural to conjecture that 
\begin{Conjecture} \label{c:K-conj}
For two-dimensional turbulent flows, the number of determining modes is bounded by $\kappa_\eta^2$.
\end{Conjecture}

In this paper we prove Conjecture~\ref{c:K-conj} for every solution of the 2D NSE, while Conjecture~\ref{c:G-conj} for non-intermittent (experimentally observable, see Section~\ref{s:Long_time_behaviour}) solutions. We show that there exists a number $\bar \lb$ (see Section~\ref{sec:Det_Wavenumber}), which is a determining wavenumber in the sense that if the difference of two solutions projected on the modes below $\bar \lb$ converges to zero, then the whole difference goes to zero in $L^2$. In particular, if two solutions on the global attractor (see Section~\ref{s:Long_time_behaviour}) coincide below $\bar \lb$, they are identical.
More precisely, we prove

\begin{Theorem}\label{thm}
Let $u(t)$ and $v(t)$ be solutions of the 2D Navier-Stokes equations (\ref{nse}).  Let $Q$ be such that $\bar\lb=\lambda_{Q}$.
If
\begin{equation} \label{eq:dm-condition}
\lim_{t \to \infty} \|u(t)_{\leq Q}- v(t)_{\leq Q}\|_{L^2} =0,
\end{equation}
then
\begin{equation}\notag
\lim_{t \to \infty} \|u(t) - v(t)\|_{L^2}=0.
\end{equation}
\end{Theorem}
Here $\bar \Lambda =\lambda_{Q} = 2^Q/L$ is the determining  wavenumber, and $u_{\leq Q}$ denotes a projection on the modes below this wavenumber (see Subsection~\ref{sec:LPD}).

Now the number of determining modes in two-dimensions is $N = \bar \Lambda^2$, which we prove satisfies the following estimate
\begin{equation} \label{e:intro-Lambda-bound}
N=\bar \lb^2 \lesssim  \kappa_\eta^2 \lesssim \kappa_0 G^{\frac{4}{d+4}}, \qquad \kappa_0=2\pi/L,
\end{equation}
with a small correction when $d$ in close to zero (but not equal to zero), which is the case of extreme intermittency not observed experimentally in turbulent flows. Here $d$ is the intermittency dimension, defined in Section~\ref{s:Long_time_behaviour}, which measures the number
of eddies at various scales.
The case $d=2$ corresponds to the regime where at each scale the eddies occupy the whole region. Even though two-dimensional turbulent flows are not expected to be highly intermittent, mathematically, $d$ could potentially be close or even equal to zero, which is called the case of extreme intermittency. Obtaining a rigorous positive lower bound on $d$ is a major open question in mathematical theory of turbulence.

While the first bound in (\ref{e:intro-Lambda-bound}) proves Conjecture~\ref{c:K-conj} modulo a corrector in a small intermittency region $d\in (0, 2\sigma)$, $\sigma \ll 1$, the second bound interpolates between two extreme cases
\begin{equation} \label{e:intro-N-G-bound}
\begin{cases}
N=\bar \lb^2 \lesssim  \kappa_0 G^{\frac{2}{3}}, \qquad &\text{for } d=2,\\
N=\bar \lb^2 \lesssim  \kappa_0 G, \qquad &\text{for } d=0,
\end{cases}
\end{equation}
where the first one reflects the optimal bound on the Hausdorff dimension of the  global attractor \eqref{e:A-dim} (with removed $\log$), while the second one recovers the estimate by John and Titi \eqref{e:A-dim} from \cite{JT}. Since experimentally and numerically two-dimensional turbulent flows are found to be non-intermittent, with a few exceptions,  estimate \eqref{e:intro-N-G-bound} validates Conjecture~\ref{c:G-conj} in the physically relevant case $d=2$.

This leads to the following
\begin{Question}
Are there solutions of the 2D NSE with nontrivial longtime behavior, which are intermittent, i.e., $d \ll 2$?
\end{Question}
 Recent experiments and numerical simulations  suggest that two-dimensional turbulent flows are not intermittent, i.e., $d=2$ \cite{BM, PT98, PJT}. There are also opposite observations \cite{CG, DR}, but even if intermittency is present, it is not expected to be strong. There is also numerical evidence that the number of determining modes does not grow as fast as the theoretical predictions \cite{OT08}. This suggests that there is a chance for Conjecture~\ref{c:G-conj} to be true for all solutions of the 2D NSE. On the other hand, the Jones-Titi bound \eqref{e:JT} and hence \eqref{e:intro-N-G-bound} might still be optimal for some pathological highly intermittent solutions of the 2D NSE, not found numerically yet, and thus  Conjecture~\ref{c:G-conj} might hold only for physically realizable turbulent solutions.

\section{Preliminaries}
\label{sec-pre}

\subsection{Notation}
\label{sec:notation}
We use the following conventional notations:
\begin{itemize}
\item $A\lesssim B$ an estimate short for $A\leq c B$ with some absolute constant $c$;
\item $A\sim B$ an estimate of the form $c_1 B\leq A\leq c_2 B$ for some absolute constants $c_1$ and $c_2$;
\item $A\lesssim_r B$ an estimate of the form $A\leq c_r B$ with some adimentional constant $c_r$ that depends on the parameter $r$;
\item $\|\cdot\|_p=\|\cdot\|_{L^p}$;
\item $(\cdot, \cdot)$ for the $L^2$-inner product.
\end{itemize}

\subsection{Littlewood-Paley decomposition}
\label{sec:LPD}
Littlewood-Paley decomposition theory are briefly recalled below, see \cite{BCD, Gr}.
For an integer $q$, $\lambda_q=2^q/L$ denotes the $q$-th shell frequeuncy. We choose a nonnegative radial function $\chi\in C_r^\infty(\R^3)$ as
\begin{equation} \label{eq:xi}
\chi(\xi):=
\begin{cases}
1, \ \ \mbox { for } |\xi|\leq\frac{3}{4}\\
0, \ \ \mbox { for } |\xi|\geq 1.
\end{cases}
\end{equation}
Denote $\varphi(\xi):=\chi(\xi/2)-\chi(\xi)$
and
\begin{equation}\notag
\varphi_q(\xi):=
\begin{cases}
\varphi(2^{-q}\xi)  \ \ \ \mbox { for } q\geq 0,\\
\chi(\xi) \ \ \ \mbox { for } q=-1.
\end{cases}
\end{equation}
Note that the sequence of $\varphi_q$ forms a dyadic partition of unity. For a tempered distribution vector field $u$ on $\T^2$, the $q$-th Littlewood-Paley projection of $u$ is defined as
\[
  u_q(x) := \Delta_q u(x) := \sum_{k\in\Z^3}\hat{u}(k)\varphi_q(k)e^{i\frac{2\pi}{L} k \cdot x},
\]
with $\hat{u}(k)$ being the $k$th Fourier coefficient of $u$. 
Thus we have $u=\sum_{q=-1}^\infty u_q$ in the distributional sense. Note that the $H^s$ norm of $u$ for $s\in\R$ satisfies
\[
  \|u\|_{H^s} \sim \left(\sum_{q=-1}^\infty\lambda_q^{2s}\|u_q\|_2^2\right)^{1/2}.
\]
We also use the short notations
\bg\notag
u_{\leq Q}:=\sum_{q=-1}^Qu_q, \quad u_{(P,Q]}:=\sum_{q=P+1}^Qu_q, \quad \tilde{u}_q := u_{q-1} + u_q + u_{q+1}.
\ed

\medskip

\subsection{Bernstein's inequality and Bony's paraproduct}
\label{sec-para}

\begin{Lemma}\label{le:bern}(Bernstein's inequality) 
Let $n$ be the spatial dimension and $1\leq s\leq r $. Then for all tempered distributions $u$ on $\mathbb T^n$, 
\bg\label{Bern}
\|u_q\|_{r}\lesssim \lambda_q^{n(\frac{1}{s}-\frac{1}{r})}\|u_q\|_{s}.
\ed
\end{Lemma}

\begin{Lemma}\label{le-bony} (Bony's paraproduct) \cite{CD-sqg} For vector fields $u$ and $v$ we have
\begin{equation}\notag
\begin{split}
\Delta_q(u\cdot\nabla v)=&\sum_{|q-p|\leq 2}\Delta_q(u_{\leq{p-2}}\cdot\nabla v_p)+
\sum_{|q-p|\leq 2}\Delta_q(u_{p}\cdot\nabla v_{\leq{p-2}})\\
&+\sum_{p\geq q-2} \Delta_q(\tilde u_p \cdot\nabla v_p).
\end{split}
\end{equation}
\end{Lemma}

\bigskip

\section{Long time behavior: Grashof and Kraichnan's numbers. Intermittency}  \label{s:Long_time_behaviour}

\subsection{Grashof's number} \label{subs:Grashof_number}
A solution $u(t)$ to the 2D NSE satisfies the enstrophy equality for the vorticity $\o= \curl u$:
\begin{equation} \label{eq:EI-sec3}
\frac{1}{2}\|\o(t)\|_2^2 = \frac{1}{2}\|\o(t_0)\|_2^2 - \nu\int_{t_0}^{t} \|\Delta u(\tau)\|_2^2\, d\tau + \int_{t_0}^{t} (\curl f,u)\, d\tau,
\end{equation}
for all $t  \geq t_0$. 

In the case of a time-independent force $f$, it is easy to see that  \eqref{eq:EI-sec3} implies
\[
\|\o(t)\|_2^2 \leq \|\o(t_0)\|_2^2 - \nu\kappa_0^2\int_{t_0}^{t} \| \o(\tau)\|_2^2\, d\tau + \frac{1}{\nu}\int_{t_0}^{t} \|f\|_{L^2}^2\, d\tau,
\]
where  $\kappa_0=2\pi\l_0 =2\pi/L$, and hence the  existence of an absorbing set
\begin{equation} \label{eq:Absorbing_ball}
B:= \{ u\in L^2(\T^3): \|\o\|_2 \leq R \}.
\end{equation}
Here the radius $R$ is such that 
\[
R > \nu \kappa_0 G,
\]
where $G$ is the adimensional Grashof number
\[
G:=\frac{\|f\|_{L^2}}{\nu^2 \kappa_0^2}.
\]
More precisely,  for any solution $u(t)$ there exists $t_0$, depending only on $\|\o(0)\|_2$, such that
\[
u(t) \in B \qquad \forall t>t_0. 
\]
Moreover, there exists a global attractor with the following structure (see \cite{FMRT}):
\[
\A = \{u(0): u(\cdot)\in L^\infty(-\infty, \infty; L^2) \text{ satisfies 2D NSE}\}.
\]
The attractor $\A \subset B$ is the  $L^2$ omega limit of $B$, and it is the  minimal $L^2$ closed attracting set. 

In the case of a time-dependent force $f=f(t)$, a relevant object describing the long-time dynamics is a pullback attractor (see \cite{CK} for a general theory applicable to the NSE). In the nonautonomous case, there exists an absorbing ball $B$ \eqref{eq:Absorbing_ball} for the vorticity $\omega$ in $L^2$, with any radius $R$ such that $R > \nu \kappa_0 G$,
just as in the autonomous case, but the Grashof number is now defined as
\begin{equation} \label{eq:Gdefin}
G:=\frac{T^{\frac12}\|f\|_{L^2_b(T)}}{\nu^{\frac 32}\kappa_0(1-e^{-\nu\kappa_0^2 T})^\frac12}= \frac{\|f\|_{L^2_b(T)}}{\nu^2 \kappa_0^{2}} \left(\frac{\mathfrak{T}}{1-e^{-\mathfrak{T}}}\right)^{\frac{1}{2}},
\end{equation}
where the averaging time is written as $T= \nu^{-1} \kappa_0^{-2} \mathfrak{T}$ for adimensional $\mathfrak{T}$ (commonly chosen to be one).
Here it is assumed that $f$ is translationally bounded in $L^2_{loc}(\mathbb{R},L^2)$ and
\[
\|f\|^2_{L^2_b(T)}:=\sup_{t\in \mathbb{R}} \frac{1}{T}\int_t^{t+T}\|f(\tau)\|_{L^2}^2 \, d\tau.
\]
The pullback attractor is defined as the minimal  closed  pullback attracting set, it is the pullback omega limit of $B$, and it has the following structure:
\[
\A(t) = \{u(t): u(\cdot)\in L^\infty(-\infty, \infty; L^2) \text{ satisfies 2D NSE}\}.
\]

Another classical way to introduce the Grashof number, as in \cite{FP, JT}, is to consider forces $f$ with
\[
F:= \limsup_{t \to \infty} \|f(t)\|_{L^2} < \infty.
\]
We cannot define the pullback attractor in this case, but can still talk about the long time behavior of solutions to the 2D NSE. Defining the Grashof number as
\[
G:=\frac{F}{\nu^2 \kappa_0^2},
\]
as in the autonomous case, we again have
\begin{equation} \label{eq:limsup_olega}
\limsup_{t \to \infty} \|\o(t)\|_{L^2} \leq  \nu \kappa_0 G,
\end{equation}
for every solution of the 2D NSE. The bound \eqref{eq:limsup_olega}, which holds for every definition of the Grashof number above, is one of the main ingredient that we need here. 

\subsection{Kraichnan's number} \label{sub:Kraichnan's_number}

Fix time $T$. For example, $T= \nu^{-1} \kappa_0^{-2}$ is a common choice, but we will keep $T$ as a free parameter.
Going back to the enstrophy equality \eqref{eq:EI-sec3} and using \eqref{eq:limsup_olega}, 
\[
\begin{split}
0\leq\|\o(t+T)\|_2^2 &\leq \limsup_{\tau \to t+}\|\o(\tau)\|_2^2 - 2\nu\int_{t}^{t+T} \|\Delta u(\tau)\|_2^2\, d\tau + 2\int_{t}^{t+T} (\curl f,u)\, d\tau\\
&\leq  (\nu\kappa_0)^2 G^2 - \nu\int_{t}^{t+T} \|\Delta u(\tau)\|_2^2\, d\tau + \frac{1}{\nu}\int_{t}^{t+T} \|f\|_{L^2}^2\, d\tau.\\
\end{split}
\]
Hence, defining the time average $\<\cdot \>$, used throughout the paper,
\[
\<F\> = \limsup_{t\to \infty} \int_t^{t+T} F(\tau) \, d\tau,
\]
and using the bound
\[
\< \|f\|_{L^2}^2 \> \leq (\nu\kappa_0)^4G^2,
\]
which clearly holds for every definition of the Grashof number in Subsection~\ref{subs:Grashof_number},
we obtain
\begin{equation} \label{eq:bound G}
\begin{split}
\<\|\Delta u\|_{L^2}^2 \>  &\leq \frac{\nu \kappa_0^2 G^2}{T} 
+ \nu^2 \kappa_0^4G^2,
\end{split}
\end{equation}

We can now connect \eqref{eq:bound G}  to  Kraichnan's number \cite{Kr}, which, taking into account (possible) intermittency effects, is defined as

\begin{equation} \label{eq:kdeps-inermit}
\kappa_\eta := \left(\frac{\eta }{\nu^3} \right)^{\frac{1}{d+4}}, \qquad  \eta := L^{-d}\nu\<\|\Delta u\|_{L^2}^2\>,
\end{equation}
where $d$ is the intermittency dimension and $\eta$ is the average enstrophy dissipation rate per unit active volume (i.e., the volume occupied by eddies). 
In order to define $d$, first note that in two dimensions, thanks to Bernstein's inequality,
\begin{equation} \label{eq:BinInt}
L^{-2}\l_q^{2} \|u_q\|_2^2 \leq  \l_q^{2} \|u_q\|_\infty^2 \leq C_{\mathrm{B}} \l_q^{2} \|u_q\|_2^4,
\end{equation}
where $C_B$ is an absolute constant, which depends on the choice of $\chi(\xi)$ in \eqref{eq:xi}. The intermittency dimension $d$ is defined as
\begin{equation} \label{eq:intermdef}
d:= \sup\left\{s\in \mathbb{R}: 0<\left<\sum_{q}\l_q^{2+s} \|u_q\|_\infty^2 \right> \leq C_{\mathrm{B}}^{2-s} L^{-s}\left<\sum_{q}\l_q^{4} \|u_q\|_2^2 \right> \right\},
\end{equation}
if the set above in nonempty, and $d=2$ in the case of the empty set, i.e, if $u$ becomes identically zero after some time.
Now \eqref{eq:BinInt} and the fact that $\<\sum_{q}\l_q^{4} \|u_q\|_2^2 \><\infty$ imply that $d \in [0,2]$ and
\[
\left<\sum_{q}\l_q^{2+d} \|u_q\|_\infty^2 \right> = C_{\mathrm{B}}^{2-d} L^{-d}\left<\sum_{q}\l_q^{4} \|u_q\|_2^2 \right>.
\]

The intermittency dimension $d$, defined in terms of  a level of saturation of Bernsten's inequality (see \cite{CD-nse,CDK,CSint} for similar definitions), measures the number
of eddies at various scales.
The case $d=2$ corresponds to the regime where at each scale the eddies occupy the whole region. Even though two-dimensional turbulent flows are not expected to be higly intermittent, mathematically, $d$ could potentially be close or equal to zero (the case of extreme intermittency). Obtaining a rigorous positive lower bound on $d$ is a major open question.

Note that the intermittency dimension $d$ and Kraichnan's number 
$\kappa_\eta$ are defined for each individual trajectory. We also define their global analogs as
\[
D:=\inf_{u } d, \qquad K_\eta:= \sup_{u} \kappa_\eta,
\]
where the infimum and supremum are taken over complete trajectories: $u \in L^\infty(\mathbb{R}; L^2)$ satisfying \eqref{nse}.

Now using the bound \eqref{eq:bound G} we can connect the Kraichnan and Grashof numbers:
\[
\kappa_\eta = \<L^{-d}\nu^{-2}\|\Delta u\|_{L^2}^2\>^{\frac{1}{d+4}}   \leq
(2 \pi)^{-\frac{d}{d+4}}\kappa_0G^{\frac{2}{d+4}}\left(\frac{1}{\nu \kappa_0^2 T} + 1\right)^{\frac{1}{d+4}}. 
\]
Also, taking the supremum over all bounded complete trajectories, with the standard choice $T= \nu^{-1} \kappa_0^{-2}$, we obtain 
\[
K_\eta \leq 2^{\frac{1}{4}}\kappa_0G^{\frac{2}{D+4}},
\]
provided $G\geq 1$.

\bigskip

\section{Definition of the determining wavenumber $\bar \Lambda$} \label{sec:Det_Wavenumber}

We first define a determining wavenumber $\lb_u(t)$ for each individual trajectory in terms of the critical norms at various scales, or local Reynolds numbers, pointwise in time, as we did for other systems in \cite{CD-sqg0, CD-nse, CDK}. Then we show that the time average $\bar \lb = \< \lb \>$ is also a determining wavenumber for the 2D NSE, which was not possible to achieve for supercritical equations. Finally, we bound  it in terms of the Kraichnan and Grashof constants.

Fix $\sigma\in(0,2]$. We define a local wavenumber $\lb^1$ in term of a critical norm with $L^2$ based for the vector field $u$
\begin{equation}\label{wave2}
\begin{split}
\lb^1_{u}(t)=\min & \left\{\lambda_q: (L\lambda_{p-q})^{\sigma}\lambda_{q}^{-1}\|\o_p\|_{L^2}<c_0\nu ,~\forall p>q \right.\\
& \left.~\text{and}~ \lambda_q^{-2}\|\nabla\o_{\leq q}\|_{L^2}<c_0\nu,~q\in \mathbb{N} \right\},
\end{split}
\end{equation}
where $c_0$ is an adimensional constant. The wavenumber $\lb^1_u$ will have an almost optimal bound when $0\leq d_u \leq 2\sigma$.

We also define a wavenumber $\lb^2$ in term of an $L^\infty$ based critical norm 
\begin{equation}\label{wave1}
\begin{split}
\lb^2_{u}(t):=\min & \left\{\lambda_q:  (L\lambda_{p-q})^{\sigma}\lambda_{q}^{-2}\|\o_p\|_{L^\infty}<c_0\nu ,~\forall p>q \right.\\
&\left. ~\text{and}~ \lambda_q^{-3}\|\nabla\o_{\leq q}\|_{L^\infty}<c_0\nu,~q\in \mathbb{N} \right\}.
\end{split}
\end{equation}
We point out that $\lb^2_u$ will enjoy an optimal bound when $d_u$ is away 
from 0. 
Note that a solution of the 2D NSEs starting from $L^2$ initial data is regular for all positive time. Hence the wavenumbers defined in (\ref{wave2}) and (\ref{wave1}) are both bounded on any time interval $[T_1, T_2]$, $T_1>0$; and, as we will see later, on $[T_1,\infty)$.

We then define
\begin{equation} \label{wave3}
\lb_u(t)=\min\{\lb^1_u(t), \lb^2_u(t)\},
\end{equation}
and 
\begin{equation} \label{wave4}
\bar \Lambda_u = \< \lb_u \>, \qquad \lb^{\mathrm{max}} = \sup_{u} \lb_u,
\end{equation}
where the supremum is taken over complete trajectories: $u \in L^\infty(\mathbb{R}; L^2)$ satisfying \eqref{nse}.

We will prove Theorem~\ref{thm}, showing that $\bar \Lambda_u$ (and hence $\lb^{\mathrm{max}}$) is a determining wavenumber in the sense that if the difference of two solutions projected on the modes below this wavenumber converges to zero, then the whole difference goes to zero in $L^2$. In particular, if two solutions on the global attractor coincide below $\lb^{\mathrm{max}}$, then they are identical.

\bigskip

\section{Proof of Theorem \ref{thm}}
\label{sec-pf1}

In this section we prove the wavenumber defined in (\ref{wave3}-\ref{wave4}) is a determining wavenumber for the 2D NSE as stated in Theorem \ref{thm}. 
Denote $\omega_1=\nabla\times u_1$ and $\omega_2=\nabla\times u_2$.
Let $v=u_1-u_2$ and $\omega=\omega_1-\omega_2$ 
which satisfies the equation
\begin{equation} \label{eq-o}
\o_t+u_1\cdot\nabla \o-\nu\Delta\o+v\cdot\nabla \o_2=0.
\end{equation}

Multiplying (\ref{eq-o}) by $\o$ and integrating yields 
\begin{equation}\label{w2}
\begin{split}
&\frac{1}{2}\frac{d}{dt}\|\o(t)\|_2^2+\nu \|\nabla\o(t)\|_2^2 \\
\leq & -\int_{\T^2}(v\cdot\nabla)\o_2 \cdot \o\, dx
-\int_{\T^2}(u_1\cdot\nabla) \o \cdot \o\, dx. 
\end{split}
\end{equation}
Note that $\int_{\T^2}(u_1\cdot\nabla) \o\cdot \o\, dx=0$. Denote $I=\int_{\T^2}(v\cdot\nabla)\o_2 \cdot \o\, dx$.
Applying the wavenumber splitting approach developed in \cite{CD-nse, CDK} we show that 
\begin{Lemma}\label{le-main}
Fix any (small) $\sigma\in(0,2)$. Let $0<\delta<1$, the wavenumber $\lb_{u_2}$ defined as in \eqref{wave3}, and $Q$ such that  $\lambda_Q=\lb_{u_2}$. Then the flux term $I$ satisfies the estimate
\begin{equation}\label{est-i1}
|I|\leq Cc_0 \nu \|\nabla \o\|_2^2+Cc_0\nu \lb_{u_2}^{4+2\delta}\sum_{p\leq Q}\l_p^{-2-2\delta} \|\o_p\|_2^2,
\end{equation}
for an absolute constant $C>0$.
\end{Lemma}
Since $\lb_{u_2}(t)=\min\{\lb^1_{u_2}(t), \lb^2_{u_2}(t)\}$, we will estimate $I$ in two separate cases according to the definition of $\lb^1_{u}$ in \eqref{wave2} and $\lb^2_{u}$ (\ref{wave1}) respectively.

\medskip

\subsection{Intermittency dimension away from zero} 
\label{sec-wave1}
Here we prove Lemma \ref{le-main} assuming $\lb^1_{u} \geq \lb^2_{u}$, and hence using the definition of wavenumber $\lb^2_{u}$ in (\ref{wave1}). We simply denote $\lb^2_{u_2}$ by $\lb$ in the following estimates without causing confusion. Let $Q$ be the integer such that $\lb=\lambda_Q=2^Q/L$.
Using Bony's paraproduct in Lemma \ref{le-bony}, $I$ can be decomposed as
\begin{equation}\notag
\begin{split}
I=&\sum_{q\geq-1}\sum_{p\geq-1}\int_{\T^2}\Delta_p((v\cdot\nabla)\o_2) \cdot \o_q\, dx\\
=&\sum_{q\geq -1}\sum_{|q-p|\leq 2}\int_{\T^2}\Delta_q(v_{\leq{p-2}}\cdot\nabla (\o_2)_p) \o_q\, dx\\
&+\sum_{q\geq -1}\sum_{|q-p|\leq 2}\int_{\T^2}\Delta_q(v_{p}\cdot\nabla (\o_2)_{\leq{p-2}})\o_q\, dx\\
&+\sum_{q\geq -1}\sum_{p\geq q-2} \int_{\T^2}\Delta_q(\tilde v_p\cdot\nabla (\o_2)_p)\o_q\, dx\\
=:&I_{1}+I_{2}+I_{3}.
\end{split}
\end{equation}
Recall that employ the definition of wavenumber in (\ref{wave1}).
Using H\"older's inequality we obtain
\[
\begin{split}
  |I_{1}| 
& \leq \sum_{q\geq -1}\sum_{|q-p|\leq 2}\int_{\T^2}|\Delta_q(v_{\leq{p-2}}\cdot\nabla (\o_2)_p) \o_q|\, dx\\
 &\lesssim  \sum_{p>Q}\sum_{|q-p|\leq 2}\|v_{(Q, p-2]}\|_{2}\lambda_p\|(\o_2)_p\|_\infty\|\o_q\|_2\\
 &+ \sum_{p>Q}\sum_{|q-p|\leq 2}\|v_{\leq Q}\|_{2}\lambda_p\|(\o_2)_p\|_\infty\|\o_q\|_2\\
 &+ \sum_{p\leq Q}\sum_{|q-p|\leq 2}\|v_{\leq p-2}\|_{2}\|\nabla (\o_2)_p\|_\infty\|\o_q\|_2\\
 &\equiv : I_{11}+I_{12}+I_{13}.
\end{split}
\]
Using the definition (\ref{wave1}) for $u_2$ and Littlewood-Paley theorem, we deduce from the H\"older's, Cauchy–Schwarz, and Jensen inequalities 
\[
\begin{split}
I_{11}  &=\sum_{p>Q}\sum_{|q-p|\leq 2}\|v_{(Q, p-2]}\|_{2}\lambda_p\|(\o_2)_p\|_\infty\|\o_q\|_2\\
&\lesssim  c_0\nu\sum_{p>Q}\sum_{|q-p|\leq 2}\lb^{2+\sigma}\lambda_p^{1-\sigma}\|\o_q\|_2\sum_{Q<p'\leq p-2}\|v_{p'}\|_{2} \\
           &\lesssim  c_0\nu\sum_{q>Q}\lambda_q\|\o_q\|_2\sum_{Q<p'\leq q}\lambda_{p'}^{2}\|v_{p'}\|_2\lambda_{p'}^{-2}\lambda_{q}^{-\sigma}\lambda_Q^{2+\sigma}\\
 &\lesssim  c_0\nu\sum_{q>Q}\lambda_q\|\o_q\|_2\sum_{Q<p'\leq q}\lambda_{p'}^{2}\|v_{p'}\|_2\lambda_{p'-q}^{\sigma}\\          
&\lesssim  c_0\nu\sum_{q>Q}\lambda_q^{2}\|\o_q\|_2^2+c_0\nu\sum_{q>Q}\left(\sum_{Q<p'\leq q}\lambda_{p'}^{2}\|v_{p'}\|_2(L\lambda_{p'-q})^{\sigma}\right)^2\\
 &\lesssim  \frac{c_0}{(1-2^{-\s})^2}\nu \|\nabla \o\|_2^2
\end{split}
\]
where we needed $\sigma>0$ to apply Jensen's inequality in the last step and used the fact $\lambda_{p'}\|v_{p'}\|_2\sim \|\omega_{p'}\|_2$.

Again using H\"older's inequality, definition (\ref{wave1}), the Cauchy–Schwarz inequality and Jensen's inequality, we have
\begin{equation}\notag
\begin{split}
 I_{12}  
 &=\sum_{p>Q}\sum_{|q-p|\leq 2}\|v_{\leq Q}\|_{2}\lambda_p\|(\o_2)_p\|_\infty\|\o_q\|_2\\
 &\lesssim  c_0\nu\sum_{p>Q}\sum_{|q-p|\leq 2}\lb^{2+\sigma}\lambda_p^{1-\sigma}\|\o_q\|_2\|v_{\leq Q}\|_{2} \\
  &\lesssim  c_0\nu\sum_{q>Q-2}\lb^{2+\sigma}\lambda_q^{1-\sigma}\|\o_q\|_2\sum_{p'\leq Q}\|v_{p'}\|_{2}\\
&\lesssim  c_0\nu\sum_{q>Q-2}\lambda_q\|\o_q\|_2\lambda_{Q-q}^{\sigma}\left(\lambda_Q^{2+\delta}\sum_{p'\leq Q}\l_{p'}^{-1-\delta}\|\o_{p'}\|_2\l_{p'-Q}^{\delta}\right)\\
&\lesssim  c_0\nu\left(\sum_{q>Q-2}\lambda_q\|\o_q\|_2(L\lambda_{Q-q})^{\sigma}\right)^2+c_0\nu\lb^{4+2\delta}\sum_{p'\leq Q}\l_{p'}^{-2-2\delta}\|\o_{p'}\|_2^2\\
&\lesssim  \frac{c_0}{(1-2^{-\s})^2}\nu \|\nabla\o\|_2^2+c_0\nu\lb^{4+2\delta}\|\nabla^{-1-\delta} \o_{\leq Q}\|_2^2
\end{split}
\end{equation}
provided $\delta>0$ (and recall $\sigma>0$).  

To estimate $I_{13}$ we first decompose it as
\[
\begin{split}
 I_{13}  
&=\sum_{q\leq Q}\sum_{|q-p|\leq 2}\|v_{\leq p-2}\|_{2}\|\nabla (\o_2)_p\|_\infty\|\o_q\|_2\\
&\leq \sum_{q\leq Q}\sum_{\substack{|q-p|\leq 2\\ p> Q}}\|v_{\leq p-2}\|_{2}\|\nabla (\o_2)_p\|_\infty\|\o_q\|_2\\
&+\sum_{q\leq Q}\sum_{\substack{|q-p|\leq 2\\ p\leq Q}}\|v_{\leq p-2}\|_{2}\|\nabla (\o_2)_p\|_\infty\|\o_q\|_2 \\
&\equiv : I_{131}+I_{132}.
\end{split}
\]
The first term is estimated as
\[
\begin{split}
I_{131}
&\leq c_0\nu\sum_{Q-1\leq q\leq Q}\sum_{\substack{|q-p|\leq 2\\ p> Q}}\lambda_Q^{2+\sigma}\lambda_p^{1-\sigma}\|v_{\leq p-2}\|_{2}\|\o_q\|_2\\
&\lesssim  c_0\nu\sum_{Q-1\leq q\leq Q}\l_q^{2}\|w_q\|_2^2+c_0\nu\lb^{4}\|v_{\leq Q}\|_{2}^2\\
&\lesssim  c_0\nu \|\nabla w\|_2^2+c_0\nu \lb^{4+2\delta}\|\nabla^{-1-\delta} \o_{\leq Q}\|_2^2
\end{split}
\]
where we used the fact that $\lambda_Q\sim\lambda_q\sim\lambda_p$ in the second step. We estimate the second term as
\[
\begin{split}
I_{132}
&\leq c_0\nu\sum_{q\leq Q}\lb^3\|\o_q\|_2\|v_{\leq q}\|_{2} \\
&\lesssim  c_0\nu\lb^{2+\delta}\sum_{q\leq Q}\l_q^{-1-\delta}\|\o_q\|_2\l_{q-Q}^{1+2\delta}\left(\lb^{2+\delta}\sum_{p'\leq q}\l_{p'}^{-1-\delta}\|\o_{p'}\|_{2}\l_{p'-q}^{\delta}\right) \\
&\lesssim  c_0\nu\lb^{4+2\delta}\sum_{q\leq Q}\l_q^{-2-2\delta}\|\o_q\|_2^2\l_{q-Q}^{2(1+2\delta)}\\
&+c_0\nu\sum_{q\leq Q}\left(\lb^{2+\delta}\sum_{p'\leq q}\l_{p'}^{-1-\delta}\|\o_{p'}\|_{2}\l_{p'-q}^{\delta}\right)^2 \\
&\lesssim  c_0\nu \lb^{4+2\delta}\|\nabla^{-1-\delta} \o_{\leq Q}\|_2^2
\end{split}
\]
provided $\delta>0$.

The term $I_2$ can be handled in an analogous way. We first decompose it 
\[
\begin{split}
  |I_{2}| &\leq \sum_{q\geq-1}\sum_{|q-p|\leq 2}\int_{\T^2}|\Delta_q(v_p\cdot\nabla (\o_2)_{\leq p-2}) \o_q| \, dx\\
&\leq  \sum_{q>Q}\sum_{|q-p|\leq 2}\|v_p\|_2\|\nabla (\o_2)_{(Q,p-2]}\|_\infty\|\o_q\|_2\\
&+\sum_{q>Q}\sum_{|q-p|\leq 2}\|v_p\|_2\|\nabla (\o_2)_{\leq Q}\|_\infty\|\o_q\|_2\\
&+\sum_{q\leq Q}\sum_{|q-p|\leq 2}\|v_p\|_2\|\nabla (\o_2)_{\leq p-2}\|_\infty\|\o_q\|_2\\
&\equiv : I_{21}+I_{22}+I_{23}
\end{split}
\]
where we adopt the convention that $(Q,p-2]$ is empty if $p-2\leq Q$. The subitems are estimated in the following
\[
\begin{split}
I_{21}
           &\leq \sum_{q> Q}\sum_{|q-p|\leq 2}\| v_p\|_2\|\o_q\|_{2}\sum_{Q<p'\leq p-2}\|\nabla (\o_2)_{p'}\|_\infty\\
           &\lesssim \sum_{q> Q}\l_q^{-1}\| \o_q\|_2^2\sum_{Q<p'\leq q}\lambda_{p'}\|(\o_2)_{p'}\|_\infty\\
           &\lesssim c_0\nu\sum_{q> Q}\l_q^{-1}\| \o_q\|_2^2\sum_{Q<p'\leq q}\lambda_{p'}^{1-\sigma}\lb^{2+\sigma}\\
           &\lesssim c_0\nu\sum_{q> Q}\l_q^{2}\| \o_q\|_2^2\sum_{Q<p'\leq q}\lambda_{p'}^{1-\sigma}\lb^{2+\sigma}\l_q^{-3}\\
           &\lesssim  c_0\nu\sum_{q> Q}\l_q^{2}\| \o_q\|_2^2
\end{split}
\]
since $\sigma>0$; and
\[
\begin{split}
I_{22}
&=\sum_{q>Q}\sum_{|q-p|\leq 2}\|v_p\|_2\|\nabla (\o_2)_{\leq Q}\|_\infty\|\o_q\|_2\\
&\lesssim  c_0\nu\sum_{q> Q}\sum_{|q-p|\leq 2}\l_q^{-1} \|\o_p\|_2\|\o_q\|_2\lb^3\\
&\lesssim  c_0\nu\sum_{q> Q}\l_q^{2} \|\o_q\|_2^2;
\end{split}
\]
\[
\begin{split}
I_{23}
&=\sum_{q\leq Q}\sum_{|q-p|\leq 2}\|v_p\|_2\|\nabla (\o_2)_{\leq p-2}\|_\infty\|\o_q\|_2\\
&\lesssim  c_0\nu\lb^3\sum_{q\leq Q}\sum_{|q-p|\leq 2} \|v_p\|_2\|\o_q\|_2\\
&\lesssim  c_0\nu\lb^3\sum_{Q-1\leq q\leq Q}\sum_{Q+1\leq p\leq Q+2} \|v_p\|_2\|\o_q\|_2\\
&+c_0\nu\lb^3\sum_{q\leq Q}\sum_{\substack{|q-p|\leq 2\\p\leq Q}} \|v_p\|_2\|\o_q\|_2\\
&\lesssim  c_0\nu\sum_{Q-1\leq q\leq Q+2}\l_q^{2} \|\o_q\|_2^2+c_0\nu\lb^{3}\sum_{q\leq Q}\l_q^{-1}\|\o_q\|_2^2\\
&\lesssim  c_0\nu \|\nabla \o\|_2^2+c_0\nu\lb^{4+2\delta}\sum_{q\leq Q}\l_q^{-2-2\delta}\|\o_q\|_2^2\l_{q-Q}^{1+2\delta}\\
&\lesssim  c_0\nu \|\nabla \o\|_2^2+c_0\nu\lb^{4+2\delta}\|\nabla^{-1-\delta}\o_{\leq Q}\|_2^2,
\end{split}
\]
where we used that $1+2\delta\geq 0$ to obtain the last step.

To estimate $I_{3}$, applying integration by parts and H\"older's inequality yields that 
\[
\begin{split}
  |I_{3}| 
  &\leq  \sum_{q\geq -1}\sum_{p\geq q-2}\int_{\T^2}|\Delta_q(\tilde v_p \otimes (\o_2)_{p}) \nabla \o_q| \, dx\\
&\leq  \sum_{p> Q}\sum_{Q<q\leq p+2}\l_q\|\tilde v_p\|_2\|\o_q\|_2\|(\o_2)_p\|_\infty\\
&+\sum_{p> Q}\sum_{q\leq Q}\l_q\|\tilde v_p\|_2\|\o_q\|_2\|(\o_2)_p\|_\infty\\
&+\sum_{p\leq Q}\sum_{q\leq p+2}\l_q\|\tilde v_p\|_2\|\o_q\|_2\|(\o_2)_p\|_\infty\\
&\equiv : I_{31}+I_{32}+I_{33}.
\end{split}
\]
By H\"older's inequality, definition (\ref{wave1}), the Cauchy-Schwarz inequality and Jensen's inequality we infer
\begin{align}\notag
\begin{split}
I_{31}
&=\sum_{p> Q}\sum_{Q<q\leq p+2}\l_q\|\tilde v_p\|_2\|\o_q\|_2\|(\o_2)_p\|_\infty\\
&\lesssim  \sum_{p> Q}\|\tilde v_p\|_2\|(\o_2)_{p}\|_\infty\sum_{Q<q\leq p+2}\lambda_q\|\o_q\|_2\\
&\lesssim  c_0\nu\sum_{p> Q}\lb^{2+\sigma}\l_p^{-\sigma}\|\tilde v_p\|_2\sum_{Q<q\leq p+2}\lambda_q\|\o_q\|_2\\
&\lesssim  c_0\nu\sum_{p> Q}\l_p\|\tilde \o_p\|_2\sum_{Q<q\leq p+2}\lambda_q\|\o_q\|_2\l_Q^{2+\sigma}\l_p^{-2-\sigma}\\
&\lesssim  c_0\nu\sum_{p> Q}\l_p\|\tilde \o_p\|_2\sum_{Q<q\leq p+2}\lambda_q\|\o_q\|_2\l_{q-p}^{2+\sigma}\\
&\lesssim  c_0\nu\sum_{p> Q}\l_p^{2}\|\o_p\|_2^2+c_0\nu\sum_{p>Q}\left(\sum_{Q<q\leq p+2}\lambda_q\|\o_q\|_2\l_{q-p}^{2+\sigma}\right)^2\\
&\lesssim  c_0\nu\sum_{p> Q}\l_p^{2}\|\o_p\|_2^2.
\end{split}
\end{align}
Similarly we have 
\begin{align}\notag
\begin{split}
I_{32}
&=\sum_{p> Q}\sum_{q\leq Q}\l_q\|\tilde v_p\|_2\|\o_q\|_2\|(\o_2)_p\|_\infty\\
&\lesssim  \sum_{p> Q}\|\tilde v_p\|_2\|(\o_2)_{p}\|_\infty\sum_{q\leq Q}\lambda_q\|\o_q\|_2\\
&\lesssim  c_0\nu\sum_{p> Q}\lb^{2+\sigma}\l_p^{-\sigma}\|\tilde v_p\|_2\sum_{q\leq Q}\lambda_q\|\o_q\|_2\\
&\lesssim  c_0\nu\sum_{p> Q}\l_p\|\tilde \o_p\|_2\left(\lb^{2+\delta}\sum_{q\leq Q}\lambda_q^{-1-\delta}\|\o_q\|_2\l_Q^{\sigma-\delta}\l_q^{2+\delta}\l_p^{-2-\sigma}\right)\\
&\lesssim  c_0\nu\sum_{p> Q}\l_p\|\tilde \o_p\|_2\left(\lb^{2+\delta}\sum_{q\leq Q}\lambda_q^{-1-\delta}\|\o_q\|_2\l_{Q-p}^{2+\sigma}\right)\\
&\lesssim  c_0\nu\sum_{p> Q-1}\l_p^{2}\|\o_p\|_2^2+c_0\nu \sum_{p> Q} \left(\lb^{2+\delta}\sum_{q\leq Q}\lambda_q^{-1-\delta}\|\o_q\|_2\l_{Q-p}^{2+\sigma}\right)^2\\
&\lesssim  c_0\nu \|\nabla \o\|_2^2+c_0\nu\lb^{4+2\delta}\|\nabla^{-1-\delta}\o_{\leq Q}\|_2^2
\end{split}
\end{align}
where we used $\sigma>0$. Turning to $I_{33}$, we first decompose it 
\[
\begin{split}
I_{33}
&=\sum_{p\leq Q}\sum_{q\leq p+2}\l_q\|\tilde v_p\|_2\|\o_q\|_2\|(\o_2)_p\|_\infty\\
&=\sum_{p\leq Q}\sum_{Q<q\leq p+2}\l_q\|\tilde v_p\|_2\|\o_q\|_2\|(\o_2)_p\|_\infty+\sum_{p\leq Q}\sum_{\substack{q\leq p+2\\q\leq Q}}\l_q\|\tilde v_p\|_2\|\o_q\|_2\|(\o_2)_p\|_\infty\\
&\equiv : I_{331}+I_{332}.
\end{split}
\]
As before, applying the definition (\ref{wave1}), H\"older's inequality, the Cauchy–Schwarz inequality and Jensen's inequality, we have
\[
\begin{split}
I_{331}
&\leq c_0\nu\sum_{p\leq Q}\lb^3\l_p^{-1}\|\tilde v_p\|_2\sum_{Q<q\leq p+2}\lambda_q\|\o_q\|_2\\
&\leq c_0\nu\sum_{Q-2<p\leq Q}\lb^3\l_p^{-1}\|\tilde v_p\|_2\sum_{Q<q\leq Q+2}\lambda_q\|\o_q\|_2\\
&\lesssim c_0\nu \lb^2\|v_{Q+1}\|_2\sum_{Q<q\leq Q+2}\lambda_q\|\o_q\|_2
+c_0\nu\sum_{Q-2<p\leq Q}\lb^3\l_p^{-1}\| v_p\|_2\sum_{Q<q\leq Q+2}\lambda_q\|\o_q\|_2\\
&\lesssim c_0\nu \sum_{Q<q\leq Q+2}\lambda_q^{2}\|\o_q\|_2^2
+c_0\nu\lb^{2+\delta}\sum_{Q-2<p\leq Q}\l_p^{-1-\delta}\| \o_p\|_2\l_{p-Q}^{-1+\delta}\sum_{Q<q\leq Q+2}\lambda_q\|\o_q\|_2\\
&\lesssim c_0\nu \sum_{Q<q\leq Q+2}\lambda_q^{2}\|\o_q\|_2^2
+c_0\nu\lb^{4+2\delta}\left(\sum_{Q-2<p\leq Q}\l_p^{-1-\delta}\| \o_p\|_2\l_{p-Q}^{-1+\delta}\right)^2\\
&\lesssim c_0\nu \|\nabla\o\|_2^2+c_0\nu\lb^{4+2\delta}\|\nabla^{-1-\delta} \o_{\leq Q}\|_2^2.
\end{split}
\]
Recalling $\tilde v_p=v_{p-1}+v_p+v_{p+1}$ and $\Lambda \|v_{Q+1}\|_2\sim \|\omega_{Q+1}\|_2$ we have 
\begin{align}\notag
\begin{split}
I_{332}
&\leq c_0\nu\sum_{p\leq Q}\lb^3\l_p^{-1}\|\tilde v_p\|_2\sum_{\substack{q\leq p+2,q\leq Q}}\lambda_q\|\o_q\|_2\\
&\lesssim c_0\nu\lb^2\|v_{Q+1}\|_2\sum_{q\leq Q}\lambda_q\|\o_q\|_2
+c_0\nu\sum_{p\leq Q}\lb^3\l_p^{-1}\|v_p\|_2\sum_{\substack{q\leq p+2, q\leq Q}}\lambda_q\|\o_q\|_2,
\end{split}
\end{align}
which we can further split using the Cauchy–Schwarz inequality,
\begin{align}\notag
\begin{split}
I_{332}&\lesssim c_0\nu\lb\|w_{Q+1}\|_2\left(\lb^{2+\delta}\sum_{q\leq Q}\l_q^{-1-\delta}\|\o_q\|_2\lambda_{q-Q}^{2+\delta }\right)\\
&+c_0\nu\lb^{4+2\delta}\sum_{ q\leq Q}\l_q^{-1-\delta}\|\o_q\|_2\l_{q-Q}^{1+2\delta}\left(\sum_{q-2\leq p\leq Q}\l_p^{-1-\delta}\| \o_p\|_2\l_{p-q}^{-1+\delta}\right)\\
&\lesssim c_0\nu\lb^{2}\|\o_{Q+1}\|_2^2+c_0\nu\lb^{4+2\delta}\left(\sum_{q\leq Q}\l_q^{-1-\delta}\|\o_q\|_2\lambda_{q-Q}^{2+\delta }\right)^2\\
&+c_0\nu\lb^{4+2\delta}\sum_{q\leq Q}\l_q^{-2-2\delta}\|\o_q\|_2^2\lambda_{q-Q}^{2(1+2\delta) }
+c_0\nu\lb^{4+2\delta}\sum_{q\leq Q}\left(\sum_{q-2\leq p\leq Q}\l_p^{-1-\delta}\| \o_p\|_2\l_{p-q}^{-1+\delta}\right)^2\\
&\lesssim c_0\nu\lb^{2}\|\o_{Q+1}\|_2^2+c_0\nu\lb^{4+2\delta}\sum_{p\leq Q}\l_p^{-2-2\delta} \|\o_p\|_2^2
\end{split}
\end{align}
where we need $-\frac12\leq \delta< 1$ to obtain the last step. 

Combining the estimates above we have for $\sigma>0$ and $0<\delta<1$ 
\begin{equation}\notag
|I|\leq Cc_0 \nu \|\nabla \o\|_2^2+Cc_0\nu \lb^{4+2\delta}\sum_{p\leq Q}\l_p^{-2-2\delta} \|\o_p\|_2^2
\end{equation}
for some absolute constant $C$.

\medskip

\subsection{Small intermittency dimension.}
\label{sec-wave2}

Now we prove Lemma \ref{le-main} assuming $\lb^1_{u} \leq \lb^2_{u}$, and hence using the definition of wavenumber $\lb^1_{u}$ in (\ref{wave2}). We denote $\lb^1_{u_2}$ by $\lb$ and let $Q$ be the integer such that $\lb=\lambda_Q=2^Q/L$.

We will estimate the flux $I=\int_{\T^2}(v\cdot\nabla)\o_2 \cdot \o\, dx$ in (\ref{w2}) in a way similar to Subsection \ref{sec-wave1}. 
The decomposition of $I$ into subitems is the same as in Subsection \ref{sec-wave1}. We only present the estimates of these subitems below. 
Using H\"older's inequality, Bernstein's inequality, definition (\ref{wave2}), the Cauchy–Schwarz  inequality and Jensen's inequality, we obtain 
\[
\begin{split}
I_{11}  &\leq \sum_{p>Q}\sum_{|q-p|\leq 2}\|v_{(Q, p-2]}\|_{\infty}\lambda_p\|(\o_2)_p\|_2\|\o_q\|_2\\
&\lesssim  c_0\nu\sum_{p>Q}\sum_{|q-p|\leq 2}\lb^{1+\sigma}\lambda_p^{1-\sigma}\|\o_q\|_2\sum_{Q<p'\leq p-2}\l_{p'}\|v_{p'}\|_{2} \\
           &\lesssim  c_0\nu\sum_{q>Q}\lambda_q\|\o_q\|_2\sum_{Q<p'\leq q}\lambda_{p'}^{2}\|v_{p'}\|_2\lambda_{p'}^{-1}\lambda_{q}^{-\sigma}\lambda_Q^{1+\sigma}\\
 &\lesssim  c_0\nu\sum_{q>Q}\lambda_q\|\o_q\|_2\sum_{Q<p'\leq q}\lambda_{p'}^{2}\|v_{p'}\|_2\lambda_{p'-q}^{\sigma}\\          
&\lesssim  c_0\nu\sum_{q>Q}\lambda_q^{2}\|\o_q\|_2^2+c_0\nu\sum_{q>Q}\left(\sum_{Q<p'\leq q}\lambda_{p'}^{2}\|v_{p'}\|_2(L\lambda_{p'-q})^{\sigma}\right)^2\\
 &\lesssim  \frac{c_0}{(1-2^{-\s})^2}\nu \|\nabla \o\|_2^2,
\end{split}
\]
provided $\sigma>0$. Similarly we have
\[
\begin{split}
 I_{12}  
 &\leq \sum_{p>Q}\sum_{|q-p|\leq 2}\|v_{\leq Q}\|_{\infty}\lambda_p\|(\o_2)_p\|_2\|\o_q\|_2\\
 &\lesssim  c_0\nu\sum_{p>Q}\sum_{|q-p|\leq 2}\lb^{1+\sigma}\lambda_p^{1-\sigma}\|\o_q\|_2\|v_{\leq Q}\|_{\infty} \\
  &\lesssim  c_0\nu\sum_{q>Q-2}\lb^{1+\sigma}\lambda_q^{1-\sigma}\|\o_q\|_2\sum_{p'\leq Q}\l_{p'}\|v_{p'}\|_{2}\\
&\lesssim  c_0\nu\sum_{q>Q-2}\lambda_q\|\o_q\|_2\lambda_{Q-q}^{\sigma}\left(\lambda_Q^{2+\delta}\sum_{p'\leq Q}\l_{p'}^{-1-\delta}\|\o_{p'}\|_2\l_{p'-Q}^{\delta+1}\right)\\
&\lesssim  c_0\nu\left(\sum_{q>Q-2}\lambda_q\|\o_q\|_2(L\lambda_{Q-q})^{\sigma}\right)^2+\frac{c_0}{1-2^{-\s}}\nu\lb^{4+2\delta}\sum_{p'\leq Q}\l_{p'}^{-2-2\delta}\|\o_{p'}\|_2^2\\
&\lesssim  \frac{c_0}{1-2^{-\s}}\nu \|\nabla\o\|_2^2+\frac{c_0}{(1-2^{-\s})}\nu\lb^{4+2\delta}\|\nabla^{-1-\delta} \o_{\leq Q}\|_2^2
\end{split}
\]
where we needed $\sigma>0$ and $\delta>-1$.  

Proceeding to $I_{13}$, we have
\[
\begin{split}
 I_{13}  
&\leq \sum_{q\leq Q}\sum_{|q-p|\leq 2}\|v_{\leq p-2}\|_{\infty}\|\nabla (\o_2)_p\|_2\|\o_q\|_2\\
&\leq \sum_{q\leq Q}\sum_{\substack{|q-p|\leq 2\\ p> Q}}\|v_{\leq p-2}\|_{\infty}\|\nabla (\o_2)_p\|_2\|\o_q\|_2\\
&+\sum_{q\leq Q}\sum_{\substack{|q-p|\leq 2\\ p\leq Q}}\|v_{\leq p-2}\|_{\infty}\|\nabla (\o_2)_p\|_2\|\o_q\|_2 \\
&\equiv : I_{131}+I_{132}.
\end{split}
\]
The estimates of $I_{131}$ and $I_{132}$ are performed as
\[
\begin{split}
I_{131}
&\leq c_0\nu\sum_{Q-1\leq q\leq Q}\sum_{\substack{|q-p|\leq 2\\ p> Q}}\lambda_Q^{1+\sigma}\lambda_p^{1-\sigma}\|v_{\leq p-2}\|_{\infty}\|\o_q\|_2\\
&\lesssim  c_0\nu\sum_{Q-1\leq q\leq Q}\l_q^{2}\|w_q\|_2^2+c_0\nu\lb^{2}\|v_{\leq Q}\|_{\infty}^2\\
&\lesssim  c_0\nu \|\nabla w\|_2^2+c_0\nu \lb^{4+2\delta}\|\nabla^{-1-\delta} \o_{\leq Q}\|_2^2
\end{split}
\]
where we used the fact $\lambda_Q\sim\lambda_q\sim\lambda_p$ again;  and
\[
\begin{split}
I_{132}
&\leq c_0\nu\sum_{q\leq Q}\lb^2\|\o_q\|_2\|v_{\leq q}\|_{\infty} \\
&\lesssim  c_0\nu\lb^{2+\delta}\sum_{q\leq Q}\l_q^{-1-\delta}\|\o_q\|_2\l_{q-Q}^{2+2\delta}\left(\lb^{2+\delta}\sum_{p'\leq q}\l_{p'}^{-1-\delta}\|\o_{p'}\|_{2}\l_{p'-q}^{1+\delta}\right) \\
&\lesssim  c_0\nu\lb^{4+2\delta}\sum_{q\leq Q}\l_q^{-2-2\delta}\|\o_q\|_2^2\l_{q-Q}^{2(2+2\delta)}\\
&+c_0\nu\sum_{q\leq Q}\left(\lb^{2+\delta}\sum_{p'\leq q}\l_{p'}^{-1-\delta}\|\o_{p'}\|_{2}\l_{p'-q}^{1+\delta}\right)^2 \\
&\lesssim  c_0\nu \lb^{4+2\delta}\|\nabla^{-1-\delta} \o_{\leq Q}\|_2^2
\end{split}
\]
provided $\delta>-1$.
The subitems of $I_2$ are estimated as
\begin{equation}\notag
\begin{split}
I_{21}
&\leq \sum_{q>Q}\sum_{|q-p|\leq 2}\|v_p\|_2\|\nabla (\o_2)_{(Q,p-2]}\|_2\|\o_q\|_\infty\\
           &\lesssim \sum_{q> Q}\sum_{|q-p|\leq 2}\| v_p\|_2\|\o_q\|_{\infty}\sum_{Q<p'\leq p-2}\|\nabla (\o_2)_{p'}\|_2\\
           &\lesssim c_0\nu\sum_{q> Q-2}\| \o_q\|_2^2\sum_{Q<p'\leq q}\lambda_{p'}^{1-\sigma}\lb^{1+\sigma}\\
           &\lesssim c_0\nu\sum_{q> Q}\l_q^{2}\| \o_q\|_2^2\sum_{Q<p'\leq q}\lambda_{p'}^{1-\sigma}\lb^{1+\sigma}\l_q^{-2}\\
           &\lesssim  c_0\nu\sum_{q> Q}\l_q^{2}\| \o_q\|_2^2
\end{split}
\end{equation}
for $\sigma>0$; and similarly 
\[
\begin{split}
I_{22}
&\leq\sum_{q>Q}\sum_{|q-p|\leq 2}\|v_p\|_2\|\nabla (\o_2)_{\leq Q}\|_2\|\o_q\|_\infty\\
&\lesssim  c_0\nu\sum_{q> Q}\sum_{|q-p|\leq 2}\l_p^{-1}\l_q \|\o_p\|_2\|\o_q\|_2\lb^2\\
&\lesssim  c_0\nu\sum_{q> Q}\l_q^{2} \|\o_q\|_2^2;
\end{split}
\]
\[
\begin{split}
I_{23}
&\leq\sum_{q\leq Q}\sum_{|q-p|\leq 2}\|v_p\|_2\|\nabla (\o_2)_{\leq p-2}\|_2\|\o_q\|_\infty\\
&\lesssim  c_0\nu\lb^2\sum_{q\leq Q}\sum_{|q-p|\leq 2} \l_q\|v_p\|_2\|\o_q\|_2\\
&\lesssim  c_0\nu\lb^2\sum_{Q-1\leq q\leq Q}\sum_{Q+1\leq p\leq Q+2} \l_q\|v_p\|_2\|\o_q\|_2\\
&+c_0\nu\lb^2\sum_{q\leq Q}\sum_{\substack{|q-p|\leq 2\\p\leq Q}}\l_q \|v_p\|_2\|\o_q\|_2\\
&\lesssim  c_0\nu\sum_{Q-1\leq q\leq Q+2}\l_q^{2} \|\o_q\|_2^2+c_0\nu\lb^{2}\sum_{q\leq Q}\|\o_q\|_2^2\\
&\lesssim  c_0\nu \|\nabla \o\|_2^2+c_0\nu\lb^{4+2\delta}\sum_{q\leq Q}\l_q^{-2-2\delta}\|\o_q\|_2^2\l_{q-Q}^{2+2\delta}\\
&\lesssim  c_0\nu \|\nabla \o\|_2^2+c_0\nu\lb^{4+2\delta}\|\nabla^{-1-\delta}\o_{\leq Q}\|_2^2
\end{split}
\]
provided $\delta\geq -1$.

In the end, applying H\"older's inequality, definition (\ref{wave2}), the Cauchy–Schwarz  inequality and Jensen's inequality, the subitems  $I_{31}$ and $I_{32}$ are estimated as follows,
\[
\begin{split}
I_{31}
&\leq  \sum_{p> Q}\|\tilde v_p\|_2\|(\o_2)_{p}\|_2\sum_{Q<q\leq p+2}\lambda_q\|\o_q\|_\infty\\
&\lesssim  c_0\nu\sum_{p> Q}\lb^{1+\sigma}\l_p^{-\sigma}\|\tilde v_p\|_2\sum_{Q<q\leq p+2}\l_q^2\|\o_q\|_2\\
&\lesssim  c_0\nu\sum_{p> Q}\l_p\|\tilde \o_p\|_2\sum_{Q<q\leq p+2}\lambda_q\|\o_q\|_2\l_Q^{1+\sigma}\l_q\l_p^{-2-\sigma}\\
&\lesssim  c_0\nu\sum_{p> Q}\l_p\|\tilde \o_p\|_2\sum_{Q<q\leq p+2}\lambda_q\|\o_q\|_2\l_{q-p}^{2+\sigma}\\
&\lesssim  c_0\nu\sum_{p> Q}\l_p^{2}\|\o_p\|_2^2+c_0\nu\sum_{p>Q}\left(\sum_{Q<q\leq p+2}\lambda_q\|\o_q\|_2\l_{q-p}^{2+\sigma}\right)^2\\
&\lesssim  c_0\nu\sum_{p> Q}\l_p^{2}\|\o_p\|_2^2,
\end{split}
\]
for $\sigma>0$. Now for $I_{32}$ we have
\[
\begin{split}
I_{32}
&\leq  \sum_{p> Q}\|\tilde v_p\|_2\|(\o_2)_{p}\|_2\sum_{q\leq Q}\lambda_q\|\o_q\|_\infty\\
&\lesssim  c_0\nu\sum_{p> Q}\lb^{1+\sigma}\l_p^{-\sigma}\|\tilde v_p\|_2\sum_{q\leq Q}\l_q^2\|\o_q\|_2.
\end{split}
\]
Using the Cauchy–Schwarz inequality, we conclude
\[
\begin{split}
I_{32} &\lesssim  c_0\nu\sum_{p> Q}\l_p\|\tilde \o_p\|_2\left(\lb^{2+\delta}\sum_{q\leq Q}\lambda_q^{-1-\delta}\|\o_q\|_2\l_Q^{\sigma-\delta-1}\l_q^{3+\delta}\l_p^{-2-\sigma}\right)\\
&\lesssim  c_0\nu\sum_{p> Q}\l_p\|\tilde \o_p\|_2\left(\lb^{2+\delta}\sum_{q\leq Q}\lambda_q^{-1-\delta}\|\o_q\|_2\l_{Q-p}^{2+\sigma}\right)\\
&\lesssim  c_0\nu\sum_{p> Q-1}\l_p^{2}\|\o_p\|_2^2+c_0\nu \sum_{p> Q} \left(\lb^{2+\delta}\sum_{q\leq Q}\lambda_q^{-1-\delta}\|\o_q\|_2\l_{Q-p}^{2+\sigma}\right)^2\\
&\lesssim  c_0\nu\sum_{p> Q-1}\l_p^{2}\|\o_p\|_2^2+c_0\nu \lb^{4+2\delta}\sum_{q\leq Q}\|\nabla^{-1-\delta}\o_q\|_2^2,
\end{split}
\]
provided $\delta\geq -3$.
We need to further split the last term $I_{33}$,
\[
\begin{split}
I_{33}
&=\sum_{p\leq Q}\sum_{q\leq p+2}\l_q\|\tilde v_p\|_2\|\o_q\|_\infty\|(\o_2)_p\|_2\\
&\lesssim  c_0\nu\sum_{p\leq Q}\lb^2\l_p^{-1}\|\tilde v_p\|_2\sum_{Q<q\leq p+2}\l_q^2\|\o_q\|_2
+ c_0\nu\sum_{p\leq Q}\lb^2\l_p^{-1}\|\tilde v_p\|_2\sum_{\substack{q\leq p+2\\q\leq Q}}\l_q^2\|\o_q\|_2\\
&\equiv : I_{331}+I_{332}
\end{split}
\]
and estimate the subitems $I_{331}$ and $I_{332}$ as
\[
\begin{split}
I_{331}
&\leq c_0\nu\sum_{Q-2<p\leq Q}\lb^2\l_p^{-1}\|\tilde v_p\|_2\sum_{Q<q\leq Q+2}\l_q^2\|\o_q\|_2\\
&\lesssim c_0\nu \|\o_{Q+1}\|_2\sum_{Q<q\leq Q+2}\l_q^2\|\o_q\|_2
+c_0\nu\sum_{Q-2<p\leq Q}\lb^2\l_p^{-1}\| v_p\|_2\sum_{Q<q\leq Q+2}\l_q^2\|\o_q\|_2\\
&\lesssim c_0\nu \sum_{Q<q\leq Q+2}\lambda_q^{2}\|\o_q\|_2^2
+c_0\nu\lb^{2+\delta}\sum_{Q-2<p\leq Q}\l_p^{-1-\delta}\| \o_p\|_2\l_{p-Q}^{-1+\delta}\sum_{Q<q\leq Q+2}\lambda_q\|\o_q\|_2\l_{q-Q}\\
&\lesssim c_0\nu \sum_{Q<q\leq Q+2}\lambda_q^{2}\|\o_q\|_2^2
+c_0\nu\lb^{4+2\delta}\left(\sum_{Q-2<p\leq Q}\l_p^{-1-\delta}\| \o_p\|_2\l_{p-Q}^{-1+\delta}\right)^2\\
&\lesssim c_0\nu \|\nabla\o\|_2^2+c_0\nu\lb^{4+2\delta}\|\nabla^{-1-\delta} \o_{\leq Q}\|_2^2,
\end{split}
\]
and
\[
\begin{split}
I_{332}
&\lesssim c_0\nu\lb\|v_{Q+1}\|_2\sum_{q\leq Q}\l_q^2\|\o_q\|_2
+c_0\nu\sum_{p\leq Q}\lb^2\l_p^{-1}\|v_p\|_2\sum_{\substack{q\leq p+2, q\leq Q}}\l_q^2\|\o_q\|_2\\
&\lesssim c_0\nu\lb\|w_{Q+1}\|_2\left(\lb^{2+\delta}\sum_{q\leq Q}\l_q^{-1-\delta}\|\o_q\|_2\lambda_{q-Q}^{3+\delta }\right)\\
&+c_0\nu\lb^{4+2\delta}\sum_{ q\leq Q}\l_q^{-1-\delta}\|\o_q\|_2\l_{q-Q}^{2+2\delta}\left(\sum_{q-2\leq p\leq Q}\l_p^{-1-\delta}\| \o_p\|_2\l_{p-q}^{-1+\delta}\right).
\end{split}
\]
Finally, thanks to the Cauchy–Schwarz inequality,
\[
\begin{split}
I_{332}&\lesssim c_0\nu\lb^{2}\|\o_{Q+1}\|_2^2+c_0\nu\lb^{4+2\delta}\left(\sum_{q\leq Q}\l_q^{-1-\delta}\|\o_q\|_2\lambda_{q-Q}^{3+\delta }\right)^2\\
&+c_0\nu\lb^{4+2\delta}\sum_{q\leq Q}\l_q^{-2-2\delta}\|\o_q\|_2^2\lambda_{q-Q}^{2(2+2\delta) }
+c_0\nu\lb^{4+2\delta}\sum_{q\leq Q}\left(\sum_{q-2\leq p\leq Q}\l_p^{-1-\delta}\| \o_p\|_2\l_{p-q}^{-1+\delta}\right)^2\\
&\lesssim c_0\nu\lb^{2}\|\o_{Q+1}\|_2^2+c_0\nu\lb^{4+2\delta}\sum_{p\leq Q}\|\nabla^{-1-\delta} \o_{p}\|_2^2,
\end{split}
\]
provided $-1\leq \delta< 1$.

Putting the estimates above together we conclude that (\ref{est-i1}) is proved for $\sigma>0$ and $0<\delta<1$.

\medskip

\subsection{Finishing the proof of Theorem \ref{thm}}
\label{sec-wave3}

First we recall a generalization of Gr\"onwall's lemma from \cite{FMTT}.

\begin{Lemma} \label{L:GenGronwall}
Let $\alpha(t)$ be a locally integrable real valued function on $(0,\infty)$, satisfying for some $0 < T < \infty$  the following conditions:
\[
\liminf_{t \to \infty} \int_t^{T+t} \phi(\tau) \, d\tau >0, \qquad \limsup_{t \to \infty} \int_t^{T+t} \phi^-(\tau) \, d\tau < \infty,
\]
where $\phi^-= \max\{ -\phi, 0\}$. Let $\psi(t)$ be a measurable real valued function on $(0,\infty)$ such that
\[
\psi(t) \to 0, \qquad \text{as} \qquad t \to \infty.
\]
Suppose $\xi(t)$ is an absolutely continuous non-negative function on $(0, \infty)$ such that
\[
\frac{d}{dt} \xi + \phi \xi \leq \psi, \qquad \text{a.e. on} \ (0,\infty).
\]
Then
\[
\xi(t) \to 0 \qquad \text{as} \qquad t \to \infty.
\]
\end{Lemma}

Now we are ready to prove the main result.
 Choosing $c_0$ small enough compared to  $C^{-1} \sim (1-2^{-\s})^{2}$, we obtain from (\ref{est-i1}) that

\begin{equation}\label{est-i2}
\begin{split}
&- \nu\|\nabla \omega\|_2^2 + I \\
&\leq 
-\frac{\nu}{2} \|\nabla \o\|_2^2+ cC_0\nu \lb^{4+2\delta}\sum_{p\leq Q}\l_p^{-2-2\delta} \|\o_p\|_2^2\\ 
&\leq -\frac{\nu}2 \bar\lb^2 \| \o_{> \bar Q}\|_2^2 + \frac{\nu}8 
\lb^{4+2\delta} \bar\lb^{-2-2\delta} \| \o_{> \bar Q}\|_2^2 + \frac{\nu}8 
\lb^{4+2\delta}  \|\o_{\leq \bar Q}\|_2^2\\
&\leq -\frac{\nu}4 \bar\lb^2 \left( \| \o\|_2^2 - 2  \| \o_{< \bar Q}\|_2^2\right) + \frac{\nu}8 
\lb^{4+2\delta} \bar\lb^{-2-2\delta} \| \o\|_2^2 + \frac{\nu}8 
\lb^{4+2\delta}  \|\o_{\leq \bar Q}\|_2^2 \\
& \leq \frac{\nu}{8} \| \o\|_2^2\bar\lb^2\left(-2+ (\lb/ \bar\lb)^{4+2\delta}\right) + \frac{\nu}8 \left( 4\bar \lb^2 + \lb^{4+2\delta} \right)  \|\o_{\leq \bar Q}\|_2^2.
\end{split}
\end{equation}
It follows from \eqref{w2} and (\ref{est-i2}) that
\begin{equation}\label{est-i3}
\frac{d}{dt}\|\omega\|_2^2 + \phi \| \omega\|_2^2 \leq \psi,
\end{equation}
where
\[
\phi(t) = \frac{\nu}8 \bar\lb^2\left(2-(\lb/ \bar\lb)^{4+2\delta} \right), \qquad
\psi(t) = \frac{\nu}8 \left( 4\bar \lb^2 + \lb^{4+2\delta} \right)  \|\o_{\leq \bar Q}\|_2^2.
\]

Since $\Lambda(t)$ is bounded, we have 
\[
\limsup_{t \to \infty} \int_t^{T+t} \phi^-(\tau) \, d\tau < \infty.
\]
Also, $\lim_{t \to \infty} w_{\leq \bar Q}(t) =0$ implies that $\lim_{t \to \infty} \psi(t) =0$. Finally, thanks to Jensen's inequality,
\begin{equation*}
\liminf_{t \to \infty} \frac{1}{T}\int_{t}^{T+t}\left( 2- (\lb/\bar \lb)^{4+2\delta} \right) \, dt\geq 2 -  \frac{1}{\bar\lb^{4+2\delta}} \limsup_{t \to \infty}\left(\frac{1}{T}\int_t^{T+t}\lb \, dt\right)^{4+2\delta}  =1 >0.
\end{equation*}
Hence
\[
\liminf_{t \to \infty} \int_t^{T+t} \phi(\tau) \, d\tau >0.
\]
Thus applying Lemma~\ref{L:GenGronwall} to (\ref{est-i3}) yields
\begin{equation}\notag
\lim_{t\to\infty}\|\omega(t)\|_2=0.
\end{equation}

\bigskip

\section{Estimates of the determining wavenumber}
\label{sec-est}

In this section we prove  estimates \eqref{e:intro-Lambda-bound} for the average wavenumber $\bar\lb$.
\subsection{Estimates in term of the Kraichnan number}
Recall that the intermittency dimension $d\in [0,2]$ is defined so that
\begin{equation} \label{eq:intermdef}
\left<\sum_{q}\l_q^{2+d} \|u_q\|_\infty^2 \right> \lesssim L^{-d}\left<\sum_{q}\l_q^{4} \|u_q\|_2^2 \right>.
\end{equation}
We first handle the case when the intermittency $d$ is away from zero and hence use definition (\ref{wave1}) of the determining wavenumber. 
Note that whenever $\lb(t) > \l_0$, we have that one of the conditions in definition of $\lb_{u}^2$ in (\ref{wave1}) is not satisfied for $q=Q-1$, i.e.,
\begin{equation}\label{eq:alt1}
2^{(p-Q+1)\sigma}\l_{Q-1}^{-2}\|\o_p\|_\infty \geq c_0\nu, \qquad \text{for some} \qquad p\geq Q,
\end{equation}
or
\begin{equation} \label{eq:alt2}
 \|\nabla \o_{\leq Q-1}\|_\infty \geq c_\s \nu\lambda_{Q-1}^3={\textstyle \frac{1}{8}}c_0 \nu\lb^3.
\end{equation}
Recall from the previous section that  $c_0 \sim (1-2^{-\s})^2$.
Thus we have
\[
(c_0\nu)^2 \lb^6 \leq  64 (\l_{p-Q}L)^{2\sigma}\lb^2 \|\o_p\|_\infty^2, \qquad \text{for some} \qquad p\geq Q, 
\]
or
\[
(c_0\nu)^2 \lb^6 \leq 64 \|\nabla \o_{\leq Q-1}\|_\infty^2,
\]
provided $\lb>\lambda_0$. Hence, adding the right hand sides, we obtain
\begin{equation} \label{eq:Lambda-main-estimate}
(c_0 \nu)^2(\lb-\l_0)^6 \lesssim  \|\nabla \o_{\leq Q-1}\|_\infty^2 + \sup_{p\geq Q}  (L\l_{p-Q})^{2\sigma}\lb^2 \|\o_p\|_\infty^2.
\end{equation}

Using \eqref{eq:Lambda-main-estimate} and Jensen's
inequality we get
\[
\begin{split}
\<\lb\>-\l_0
& \lesssim \<(\lb - \l_0)^{d+4}\>^{\frac{1}{d+4}}\\
& \lesssim \left< \frac{\lb^{d-2}}{(c_0 \nu)^2}\left( \|\nabla \o_{\leq Q-1}\|_\infty^2 + \sup_{q\geq Q}  (L\l_{q-Q})^{2\sigma} \lb^2\|\o_q\|_\infty^2\right)\right>^{\frac{1}{d+4}}\\
&\lesssim \left\< \frac{1}{(c_0\nu)^2}\left(\sum_{q\leq Q-1} \l_q^{\frac{d}{2}}\|\o_q\|_\infty (L\lambda_{Q-q})^{\frac{d-2}{2}} \right)^2
+ \frac{\lb^{d}}{\nu^2}\sup_{q\geq Q}  (L\l_{q-Q})^{2\sigma} \|\o_q\|_\infty^2
\right\>^{\frac{1}{d+4}}\\
&\lesssim \left\<\frac{1}{(c_0\nu)^2} \sum_{q\leq Q-1} \l_q^{d}\|\o_q\|_\infty^2 + \frac{1}{\nu^2}\sup_{q\geq Q}  (L\l_{q-Q})^{2\sigma-d} \l_q^{d}\|\omega_q\|_\infty^2
 \right\>^{\frac{1}{d+4}}\\
&\lesssim \left\<\frac{c_0^{-2}+1}{\nu^2} \sum_{q} \l_q^{d}\|\o_q\|_\infty^2 
\right\>^{\frac{1}{d+4}},
\end{split}
\]
where $d\in[2\sigma, 2]$ was used at the last step. Since $\|\o_q\|_\infty \lesssim \lambda_q \|u_q\|_\infty$ by Bernstein's inequality,
it then follows from (\ref{eq:intermdef})
\[
\<\lb\>-\l_0 \lesssim \left\<\frac{\nu }{\nu^3} \sum_{q} \l_q^{2+d}\|u_q\|_\infty^2 
\right\>^{\frac{1}{d+4}}\lesssim \left\<\frac{L^{-d}\nu}{\nu^3} \sum_{q} \l_q^{4}\|u_q\|_2^2 
\right\>^{\frac{1}{d+4}} = \kappa_\eta.
\]

\noindent
{\bf Extreme intermittency.} Now we consider the case  $d\in[0, 2\sigma]$, not expected in turbulent flows, but still plausible mathematically.
For $\lb(t) > \l_0$, one of the conditions in definition of $\lb_{u}^1$ in (\ref{wave2}) is not satisfied for $q=Q-1$, i.e.,
\begin{equation}\label{eq:alt3}
2^{(p-Q+1)\sigma}\l_{Q-1}^{-1}\|\o_p\|_2 \geq c_0\nu, \qquad \text{for some} \qquad p\geq Q,
\end{equation}
or
\begin{equation} \label{eq:alt4}
 \|\nabla \o_{\leq Q-1}\|_2 \geq c_0 \nu\lambda_{Q-1}^2={\textstyle \frac{1}{4}}c_\s \nu\lb^2.
\end{equation}
It follows 
\[
(c_0\nu)^2 \lb^4 \leq  16 (\l_{p-Q}L)^{2\sigma}\lb^2 \|\o_p\|_2^2, \qquad \text{for some} \qquad p\geq Q, 
\]
or
\[
(c_0\nu)^2 \lb^4 \leq 16 \|\nabla \o_{\leq Q-1}\|_2^2.
\]
Combining the two estimates above yields
\begin{equation} \label{eq:Lambda-main-estimate2}
(c_0 \nu)^2(\lb-\l_0)^4 \lesssim  \|\nabla \o_{\leq Q-1}\|_2^2 + \sup_{p\geq Q}  (L\l_{p-Q})^{2\sigma}\lb^2 \|\o_p\|_2^2.
\end{equation}

We infer from \eqref{eq:Lambda-main-estimate2} and Jensen's
inequality 
\[
\begin{split}
\<\lb\>-\l_0
& \lesssim \<(\lb - \l_0)^{4}\>^{\frac{1}{4}}\\
& \lesssim \left< \frac{1}{(c_0 \nu)^2}\left( \|\nabla \o_{\leq Q-1}\|_2^2 + \sup_{q\geq Q}  (L\l_{q-Q})^{2\sigma} \lb^2\|\o_q\|_2^2\right)\right>^{\frac{1}{4}}\\
&\lesssim \left\< \frac{1 }{\nu^2}\sum_{q\leq Q-1} \l_q^{2}\|\o_q\|_2^2 
+ \frac{1}{\nu^2}\sup_{q\geq Q}  (L\l_{q-Q})^{2\sigma-2} \lambda_q^2 \|\o_q\|_2^2
\right\>^{\frac{1}{4}}\\
&\lesssim \left\<\frac{1}{\nu^2} \sum_{q} \l_q^{2}\|\o_q\|_2^2 
\right\>^{\frac{1}{4}}
\end{split}
\]
since $\sigma<1$. 
Thanks to $\|\o_q\|_2 \lesssim \lambda_q \|u_q\|_2$ by Bernstein's inequality, we thus have for $d\leq 2\sigma$,
\[
\<\lb\>-\l_0 \lesssim  \left\<\frac{1}{\nu^2} \sum_{q} \l_q^{4}\|u_q\|_2^2 
\right\>^{\frac{1}{4}} = \kappa_\eta  \left(\frac{\kappa_\eta}{\kappa_0}\right)^{\frac{d}{4d+16}} \leq \kappa_\eta  \left(\frac{\kappa_\eta}{\kappa_0}\right)^{\frac{d}{16}}.
\]
where we assumed that $\kappa_\eta \geq \kappa_0$ at the last step.
To conclude, when $\kappa_\eta \geq \kappa_0$, we have
\begin{equation} \label{eq:K-estimates}
\begin{cases}
\bar \lb \lesssim  \kappa_\eta,\qquad &\text{for } d\geq 2\sigma,\\
\bar \lb \lesssim  \kappa_\eta \left(\frac{\kappa_\eta}{\kappa_0}\right)^{\frac{d}{16}}, \qquad &\text{for } d <2\sigma .
\end{cases}
\end{equation}

In the above estimate we assumed that $\sigma$ was a fixed small positive number. If we are interested in the extreme intermittency regime $d \in(0, 2\sigma)$, not observed experimentally, then recall that Theorem~\ref{le-main} holds for any $\sigma \in (0,2)$ in the definition of $\Lambda^2_u$ as long as $c_0 \lesssim (1-2^{-\sigma})^2$. Then given the intermittency dimension $d>0$, the optimal choice for $\sigma$ is $\sigma = d/2$, and we arrive at a more precise estimate
\[
\bar \lb \lesssim  \kappa_\eta \min\left\{(1-4^{-d})^{-\frac{4}{d+4}}, \left(\kappa_0^{-1}\kappa_\eta\right)^{\frac{d}{4d+16}} \right\}.
\]

\subsection{Estimates in term of the Grashof number}

Recall from Subsection~\ref{sub:Kraichnan's_number} that 
\[
\kappa_\eta \lesssim \kappa_0 G^{\frac{2}{d+4}},
\]
provided $G\geq1$. Thus
\begin{equation} \label{eq:G-estimates}
\begin{cases}
\bar \lb \lesssim  \kappa_0 G^{\frac{2}{d+4}}, \qquad &\text{for } d\geq 2\sigma,\\
\bar \lb \lesssim  \kappa_0 G^{\frac{2}{d+4}} G^{\frac{d}{2(d+4)^2}}, \qquad &\text{for } d <2\sigma ,
\end{cases}
\end{equation}

Note that the extreme intermittency regime $d=0$ is the worst case scenario for the estimates \eqref{eq:G-estimates}. Indeed, \eqref{eq:G-estimates} implies 
\[
\bar \lb \lesssim  \kappa_0 G^{\frac{1}{2}},
\]
i.e., the number of determining modes ${\bar \lb}^2$ is always bounded by $G$, which recovers the result by Jones and Titi \cite{JT}.

In the physically relevant case $d=2$, usually observed in experiments, the estimates \eqref{eq:G-estimates} become
\[
\bar \lb\lesssim  \kappa_0 G^{\frac13}
\]
which is consistent with the classical result \cite{CFT} on the  finite dimensionality of the  global attractor $\mathcal A$
 \[
 \mathrm{dim_H}\mathcal A \lesssim G^\frac{2}{3}(1+\log G),
 \]
 with the logarithmic factor of $G$ removed.

\bigskip

\section*{Acknowledgement}
A. Cheskidov was partially supported by the NSF grant DMS-1909849.
M. Dai was partially supported by the NSF grants DMS-1815069 and DMS-2009422, and the von Neumann Fellowship.  The authors are grateful to IAS for its hospitality in 2021-2022 during which time the work was completed.

\bigskip

\bibliographystyle{alpha}

\end{document}